\documentstyle[12pt]{article}

\begin{document}

\newcommand{\pf}{

\smallskip

\noindent {\it Proof : }}

\newcommand{\ti}[1]{\mbox{$\tilde{#1}$}}

\newcommand{\wti}[1]{\mbox{$\widetilde{#1}$}}

\newcommand{\T}{\mbox{$\tilde{T}$}}

\newcommand{\U}{\mbox{$\tilde{U}$}}

\newcommand{\nul}{\mbox{$\tilde{0}$}}

\newcommand{\al}{\mbox{$\alpha$}}

\newcommand{\be}{\mbox{$\beta$}}

\newcommand{\dd}{\mbox{$\delta$}}

\newcommand{\Zo}{\mbox{$Z^{\perp}$}}

\newcommand{\Yo}{\mbox{$Y^{\perp}$}}

\newcommand{\la}{\mbox{$\lambda$}}

\newcommand{\e}{\mbox{$\epsilon$}}

\newcommand{\norm}[1]{\mbox{$\|#1\|$}}

\newcommand{\cvf}{\mbox{$\stackrel{w}{\rightarrow}$}}

\newcommand {\ram}{[\omega]^\omega}

\newcommand {\ca} {2^\omega}

\newcommand{\cantor}{2^{<\omega}}

\newcommand {\N}{\mathbb N}
\newcommand {\poset}{{\mathcal P}}
\newcommand {\Q}{\mathbb Q}

\newcommand {\R}{\mathbb R}

\newcommand {\Z}{\mathbb Z}

\newcommand {\C}{\mathbb C}

\newcommand{\om}{\omega}

\newcommand{\eps}{\epsilon}

\newcommand{\compl}{\complement}

\newcommand {\iso} {\cong}

\newcommand{\isom}{\simeq}

\newcommand{\ext}{\preccurlyeq}

\newcommand {\tom} {\emptyset}

\newcommand{\emb}{\sqsubseteq}

\newcommand{\inj}{\hookrightarrow}

\newcommand{\begr}{\upharpoonright}

\newcommand{\surj}{\twoheadrightarrow}

\newcommand{\bij}{\longleftrightarrow}

\newcommand{\hviss}{\longleftrightarrow}

\newcommand{\hvis}{\Longleftarrow}

\newcommand{\saa}{\Longrightarrow}

\newcommand{\til}{\longrightarrow}

\newcommand{\equi}{\Longleftrightarrow}

\newcommand{\Lim}[1]{\mathop{\longrightarrow}\limits_{#1}}

\newcommand {\Del}{ \; \Big| \;}

\newcommand {\del}{ \; \big| \;}

\newcommand {\Mgdv}{\Big\{}

\newcommand {\mgdv}{\big\{}

\newcommand {\Mgdh}{\Big\}}

\newcommand {\mgdh}{\big\}}

\newcommand {\Intv}{\Big[}

\newcommand {\intv}{\big[}

\newcommand {\Inth}{\Big]}

\newcommand {\inth}{\big]}

\newcommand {\For}{\Bigcup}

\newcommand {\for}{\bigcup}

\newcommand {\Snit}{\Bigcap}

\newcommand {\snit}{\bigcap}

\newcommand {\og}{\; \land \;}

\newcommand {\eller}{\; \vee\;}

\newcommand{\ikke}{\lnot}

\newcommand {\go} {\mathfrak}

\newcommand {\ku} {\mathcal}

\newcommand {\un} {\underline}

\newcommand {\ex} {\exists}

\renewcommand {\a} {\forall}

\newcommand{\fed}{\mathbf}

\font\tenBbb=msbm10  \font\sevenBbb=msbm7  \font\fiveBbb=msbm5

\newfam\Bbbfam

\textfont\Bbbfam=\tenBbb \scriptfont\Bbbfam=\sevenBbb

\scriptscriptfont\Bbbfam=\fiveBbb

\def\Bbb{\fam\Bbbfam\tenBbb}

\def\R{{\Bbb R}}

\def\C{{\Bbb C}}

\def\N{{\Bbb N}}

\def\zed{{\Bbb Z}}

\def\Q{{\Bbb Q}}

\newcommand{\pff}{$\hfill  \Box$

\smallskip }

\newtheorem{defi}{Definition}

\newtheorem{prop}[defi]{Proposition}

\newtheorem{lemm}[defi]{Lemma}

\newtheorem{theo}[defi]{Theorem}

\newtheorem{coro}[defi]{Corollary}

\title{ERGODIC BANACH SPACES}

\author{V. Ferenczi and C. Rosendal}

\date{ }

\maketitle

\begin{abstract}{\em We show that any Banach space contains
 a continuum of non isomorphic subspaces or a minimal subspace. We define
an ergodic
Banach space $X$ as a space such that $E_0$
Borel reduces to isomorphism on the set of subspaces of $X$, and show that
every Banach space is
either ergodic or contains a subspace with an unconditional basis which is
complementably universal for
the family of its block-subspaces.
We also use our methods to get
uniformity results;
for example, in combination with a result of B. Maurey, V. Milman and N.
Tomczak-Jaegermann, we show that an
 unconditional basis of a Banach space, of which
every block-subspace is complemented, must be
asymptotically 
$c_0$ or
$l_p$.}
\end{abstract}

 The following question was asked the authors by G. Godefroy: how many
non isomorphic subspaces must a given Banach space contain?
By the results of W.T. Gowers, W.T. Gowers - B. Maurey and
 R. Komorowski- N. Tomczak solving the homogeneous space problem,
if $X$ is not isomorphic to $l_2$ then it must contain at least
two non-isomorphic subspaces. Except $l_2$, no examples of spaces with
only finitely, or even countably many, isomorphism classes of subspaces
 are known,
so we may ask what the possible number of non-isomorphic subspaces of a given
Banach space is, supposing it being non isomorphic to $l_2$.
This question may also be asked in the setting of the classification
of analytic equivalence relations up to Borel reducibility. If $X$ is not
isomorphic
to $l_2$, when can we classify the relation of isomorphism on subspaces of $X$?

A topological space $X$ is said to be Polish if it is separable and its
topology can be generated
 by a complete metric. Its Borel subsets are those belonging to the
smallest $\sigma$-algebra containing the
open sets. A subset is analytic if it is the continuous direct image of a
Polish space or equivalently of a
Borel set in a Polish space. All uncountable Polish spaces turn out to be
Borel isomorphic, i.e., isomorphic
by a function that is Borel bimeasurable.
 A $C$-measurable set is one belonging to the smallest $\sigma$-algebra
containing the open sets and closed under the Souslin operation, in
particular all analytic sets are $C$-measurable.

All $C$-measurable sets are universally measurable, i.e., measurable with
respect
to any
$\sigma$-finite Borel
 measure on the space. Furthermore, they have the Baire property, i.e.,
can be written on the form
$A=U\triangle M$, where $U$ is open and $M$ is meagre and are completely
Ramsey. In fact they satisfy
almost any regularity property
satisfied by Borel sets.

Most results contained in this article are centered around the notion of
Borel reducibility.
 This notion turns out to be extremely useful as a mean of measuring
complexity in analysis. It also gives
another refined view of cardinality, in that it provides us with a notion
of the number of classes of an
equivalence relation before everything gets muddled up by the wellorderings
provided by the axiom of choice.

\begin{defi}Suppose that $E$ and $F$ are analytic equivalence
  relations on Polish spaces $X$ and $Y$ respectively, then we write
$E\leq_BF$ iff there is a Borel function
$f:X\til Y$, such that $xEy\hviss f(x)Ff(y)$. Moreover, we denote by
$E\sim_BF$ the fact that the relations
are Borel bireducible, i.e., $E\leq_BF$ and $F\leq_BE$. \end{defi}

Then $E\leq_BF$  means that there is an injection from $X/E$ into $Y/F$
admitting a Borel lifting.
 Intuitively, this says that the objects in $X$ are simpler to classify
with respect to $E$ than the objects
in $Y$ with respect to $F$. Or again that $Y$ objects modulo $F$ provide
complete invariants for $X$ objects
with respect to $E$-equivalence, and furthermore, these invariants can be
calculated in a Borel manner from
the initial objects.

We call an equivalence relation $E$  on a Polish space $X$ {\em smooth} if
it Borel
 reduces to the identity relation on $\R$, or in fact to the identity
relation on any uncountable Polish
space. This is easily seen to be  equivalent to admitting a countable
separating family $(A_n)$ of Borel
sets, i.e., such that for any $x,y\in X$ we have $xEy\equi \a n\;( x\in
A_n\hviss y\in A_n)$.

A Borel probability measure $\mu$ on $X$ is called {\em $E$-ergodic} if for
any $\mu$-measurable $A\subset X$
that is $E$-invariant, i.e., $x\in A \og xEy\til y\in A$, either $\mu
(A)=0$ or $\mu (A)=1$. We call $\mu$
{\em $E$-non atomic} if every equivalence class has measure $0$.

Suppose $\mu$ was $E$-ergodic and $(A_n)$ a separating family for $E$. Then
by ergodicity and the fact that the $A_n$ are invariant either $A_n$ or $
A^c_n$ has measure $1$, so $\snit\mgdv A_n\del \mu (A_n)=1\mgdh\cap \snit
\mgdv  A ^c_n\del \mu (A_n)=0\mgdh$ is an $E$ class of full measure and
$\mu$ is
atomic. So a smooth equivalence relation cannot carry an ergodic, non
atomic probability measure.

The minimal non smooth Borel equivalence relation is the relation of
eventual agreement
of infinite binary sequences, $E_0$. This is defined on $\ca=\{0,1\}^\N$
by
$$xE_0y\hviss \ex n\; \a m\geq n\; x_m=y_m$$
To see that $E_0$ is non smooth just notice that the usual coinflipping
measure on$\ca$ is $E_0$ atomic and ergodic by the zero-one law.
Furthermore, any perfect set of almost disjoint infinite subsets of $\N$
shows that
$E_0$ has a perfect set of classes.

If $E$ is an equivalence relation on a set $X$
and $A\subset X$ then we call $A$ a {\em transversal for $E$ on $X$} if it
intersects every $E$-equivalence
class in exactly one point.
We notice that if $E$ is an equivalence relation and $A$ a transversal for
$E$, both of them analytic, then $E$ is  smooth.

Our general reference for descriptive set theory and Ramsey theory is
\cite{kec}, we adopt his notation wholesale.

\

 It is  natural to try to distinguish some
class of Banach spaces by a condition
 on the number of non isomorphic subspaces.  A step up from homogeneity
would be when  the subspaces would at
least admit some  classification in terms of real numbers, i.e., something
ressembling type or entropy. This
would say that in some sense the space could not be too wild and one would
expect such a space to have
more regularity
properties than those of a more generic space, in particular than those of
a hereditarily indecomposable
space.

With the discoveries in Banach space theory in the seventies and nineties
the hope for a general
 classification of the subspaces of a Banach space faded and even the hope
of finding some nice class of
subspaces contained in every Banach space also turned out to be doomed to
failure, by the examples of Tsirelson
and Gowers-Maurey. The strongest positive result was Gowers' dichotomy
saying that either a Banach space
contains a   hereditarily indecomposable subspace  or 
subspace with an unconditional
basis, that is either a very
rigid space (with few isomorphisms and projections) or a somewhat nice
space (with many isomorphisms and
projections).

We  isolate another class of separable Banach spaces, namely those on which
the isomorphism relation between subspaces does not reduce $E_0$, the non
ergodic ones, in particular this class
includes those admitting classification by real numbers, and show that if a
space belongs to this class, then
it must satisfy some useful regularity properties. Let ${\cal B}_X$ be the
space of closed linear subspaces of
a Banach space $X$, equipped with
 its Effros-Borel structure (see \cite{kec} or \cite{FR}). We note that
isomorphism is analytic on ${\cal B}_X^2$.
Let us define a Banach space $X$ to be {\em ergodic} if the relation $E_0$
Borel reduces to isomorphism on subspaces of $X$.
In \cite {FR}, \cite{R}, the authors studied spaces generated
by subsequences of a space $X$ with a basis: for $X$  a Banach space
with an unconditional basis, either
$X$ is ergodic  
or $X$ is isomorphic to its hyperplanes, to its square, and more
generally to any direct sum $X
\oplus Y$ where
$Y$ is generated by a subsequence of the basis,  and satisfies other regularity
properties.

Note that it is easily checked that Gowers' construction of a space with a
basis, such that no disjointly supported
subspaces are isomorphic (\cite{G0},\cite{GM2}), provides an example of a
space for which the complexity of
isomorphism on subspaces generated by subsequences is exactly $E_0$.

\

In the main part of this article,
we shall consider  a Banach space with a basis, and restrict our attention
to subspaces generated by block-bases.
As long as we consider only block-subspaces, there are more examples
of spaces with low complexity, for example
$l_p$, $1 \leq p$ or $c_0$ has only one class of isomorphism for
block-subspaces.
After noting a few facts about the number of non-isomorphic subspaces of a
Banach space, that come as consequences of Gowers' dichotomy theorem (Lemma
\ref{uniformity} to
Theorem \ref{minimality}),
we prove that block-subspaces in a non-ergodic Banach space satisfy
regularity properties
(Theorem \ref{theo7}, Corollary \ref{coro8}, Theorem \ref{theo7.5}).
We then show how our methods yield uniformity results and in particular find
new  characterizations of the classical spaces $c_0$ and $l_p$ (Corollary
\ref{asymptotic}).
Finally we show how to generalize our results to subspaces with a finite
dimensional decomposition on the basis (Theorem \ref{theo19}, Proposition
\ref{prop20}).

\

Let us recall a definition from H. Rosenthal: we say that a space $X$ is
{\em minimal} if $X$ embeds in any of its subspaces. Minimality is hereditary.
In the context of block-subspaces, there are two natural definitions:
 we define a space $X$ with a basis to be
{\em block-minimal} if every  block-subspace of $X$ has a further
block-subspace isomorphic to $X$; it is
{\em equivalence block-minimal}
if every  block-subspace of $X$ has a further
block-subspace equivalent to $X$.
The second property is hereditary, but the first one is not, so we also
define
a {\em hereditarily block-minimal} space
 as a space $X$ with a basis such that any of its block-subspaces is
block-minimal.

Let $X$ be a Banach space with a basis $\{e_i\}$. 
If
$y=(y_n)_{n
\in
\N}$ is a block-sequence of
$X$, we denote by $Y=[y_n]_{n \in
\N}=[y]$ the closed linear span of $y$. For two finite or
 infinite block-bases $z$ and $y$ of $\{e_i\}$, write
$z\leq y$ if $z$ is a blocking of $y$ (and write $Z \leq Y$ for the
corresponding subspaces). If $y=(y_i)_{i
\in \N}$, $z=(z_i)_{i \in \N}$ and
$N\in\N$, write $z\leq^* y$
iff there is an $N$ such that $(z_i)_{i \geq N} \leq y$ (and write $Z
\leq^* Y$ for the corresponding subspaces). If
$s=(s_1,\ldots,s_n)$ and $t=(t_1,\ldots,t_k)$ are two finite block-bases, i.e.,
$supp(s_i)<supp(s_{i+1})$ and  $supp(t_i)<supp(t_{i+1})$, then we write $s\preceq t$ iff $s$ is an
initial segment of $t$, i.e.,
$n\leq k$ and $s_i=t_i$ for $i\leq n$. In that case we write $t\setminus s$ for $(t_{n+1},\ldots,t_k)$.
If $s$ is a finite block-basis and $y$ is a finite or infinite block-basis
supported after $s$, denote by $s^{\frown}y$  the concatenation of 
$s$ and $y$. 

We denote by $bb(X)$ the set of normalized block-bases on $X$.
This set can be equipped with the product topology of the norm topology on $X$,
 in which case it becomes a Polish space that we denote by $bb_N(X)$.

Sometimes we want to work with blocks with rational coordinates,
though we no longer
 can demand these to be normalized (by rational, we shall always mean
an element of
$\Q+i\Q$ in the case of a complex Banach space). We identify the set of
such blocks with the set
$\Q^{<\om}_*$ of finite, not identically zero, sequences of rationals. 
We shall denote by $(\Q^{<\om}_*)^\om$ the set of (not
necessarily successive)  infinite sequences of rational blocks. Again
when needed we will give 
$\Q^{<\om}_*$ the discrete topology and 
$(\Q^{<\om}_*)^\om$ the product topology. The set of rational
block-bases may
be seen as a subset of
$(\Q^{<\om}_*)^\om$ and is denoted by $bb_\Q$. The set of finite
rational block-bases is then denoted by $fbb_\Q$.

Finally for the topology that interests us the most: let  $\bf {Q}$ be the set 
of normalized blocks of the basis that are a multiple of some block with
rational coordinates; we
denote by $bb_d(X)$ the set of block-bases of vectors in $\bf {Q}$, equipped with the product topology
of the discrete topology on $\bf {Q}$.  As $\bf {Q}$ is countable, this topology is Polish and epsilon matters
may be forgotten until the applications; when we deal with isomorphism classes, they are not relevant since a
small enough perturbation preserves the class. Note also that the canonical embedding of
$bb_d(X)$ into
${\cal B}_X$ is Borel, and this allows us to forget about the Effros-Borel structure
when checking ergodicity. Unless specified otherwise,
 from now on we work with this topology.

We first prove a Lemma about uniformity for these properties. For $C \geq 1$,
we say that a space $X$ with a basis is {\em $C$ block-minimal (resp. $C$
equivalence block-minimal)} if any block-subspace of $X$ has a further
block-subspace
which is $C$-isomorphic (resp. $C$-equivalent) to $X$.

\begin{lemm} \label{uniformity}(i)  Let $X$ be a Banach space with a basis
and assume
$X$ is (equivalence) block-minimal. Then there exists $C \geq 1$ such that
$X$ is $C$ (equivalence) block-minimal.

(ii) Suppose $\{e_n\}$ is
a basis in a Banach space, such that any subsequence of $\{e_n\}$ has a
block-sequence equivalent to
$\{e_n\}$. Then there is a subsequence $\{f_n\}$ of $\{e_n\}$ and a
constant $C \geq 1$, such that
any subsequence of $\{f_n\}$ has a block-sequence $C$-equivalent to
$\{e_n\}$.
\end{lemm}

\pf We will only prove (i) as the proof of (ii) is similar.
 Let
for $n \in \N$ $c(n)$ denote a constant such that
any $n$-codimensional subspaces of any Banach space are $c(n)$-isomorphic
(Lemma 3
in \cite{FR}). Let $X$ be block-minimal.
We want to construct by induction a decreasing sequence of block-subspaces
$X_n, n \geq 1$
 and successive
block-vectors $x_n$ such that the first vectors of $X_n$ are
$x_1,\ldots,x_{n-1}$ and such that
no block-subspace of $X_n$ is $n$ isomorphic to $X$.
Assume we may carry out the induction: then for all $n \in \N$,
no block-subspace of $[x_n]_{n \in \N}$ is $n$-isomorphic
to $X$, and this contradicts the block-minimality of $X$.
So the induction must stop at some $n$, meaning that
every block-subspace of $X_n$ whose first vectors are
$x_1,\ldots,x_{n}$ has a further block-subspace $n$ isomorphic
to $X$.
Then by definition of $c(n)$, every block-subspace of $X_n$ has a further
block-subspace $nc(n)$ isomorphic to $X$.
By block-minimality we may assume that $X_n$ is $K$-isomorphic to $X$ for some
$K$.
Take now any block-subspace $Y$ of $X$, it is $K$-isomorphic to a subspace
of $X_n$;
by standard perturbation arguments, we may find a block-subspace of $Y$ which is
$2K$-equivalent
 to
a block-subspace of $X_n$, and by the above, an even further block-subspace
  $2K$ equivalent to a $nc(n)$-isomorphic
copy of $X$; so finally $Y$ has a block-subspace $2Knc(n)$ isomorphic to
$X$ and so,
$X$ is $2Knc(n)$ block minimal.

We may use the same proof for equivalence block-minimality, using instead
of $c(n)$ a
constant $d(n)=(1+(n+1)c)^2$, such that any two normalized block-sequences
differing by only
the $n$ first vectors are $d(n)$-equivalent ($c$ stands for the constant of
the basis).
 \pff

\

 Let us recall a version of the Gowers' game $G_{A,Y}$ shown to be
equivalent to Gowers' original game by Bagaria and Lopez-Abad:
Player I plays in the $k$'th move a normalized block vector $y_k$ of $Y$
such that $y_{k-1}<y_k$ and Player II responds by either doing nothing or
playing a normalized block vector $x\in [y_{i+1},\cdots,y_k]$  if $i$ was
the last move where she played a vector. Player II wins the game if in the
end she has produced an infinite sequence $(x_k)_{k \in \N}$ which is a
block sequence in $A$. If Player II has a winning
strategy for $G_{A,Y}$ we say that she  has a
winning strategy for Gowers' game in
$Y$ for producing block-sequences in $A$.
Gowers proved that if $A$ is analytic in $bb_N(X)$,  such that any
normalized block sequence contains a further normalized block sequence
in $A$, then II has a winning strategy in
some $Y$ to produce a block-sequence arbitrarily close to a block-sequence
in $A$.

As an application of Gowers' theorem one can mention that if $X$ is (resp.
equivalence) block-minimal, then there is a constant $C$, such that for
every block
subspace $Y\leq X$, Player II has a winning strategy for Gowers' game in
some $Z \leq Y$
for producing block-sequences
 spanning a space $C$-isomorphic to $X$
 (resp. $C$-equivalent to the basis of $X$).

\

We recall that a space with a basis is said to be {\em quasi-minimal} if
any two
block-subspaces
have further isomorphic block-subspaces.
On the contrary, two spaces are said to be {\em totally incomparable} if
no subspace of the first one is isomorphic to a subspace of the second.
Using his dichotomy theorem,
Gowers proved the following result about Banach spaces.

\begin{theo} \label{gowers}
(Gowers' "trichotomy") Let $X$ be a Banach space.
Then
$X$ either contains

- a hereditarily indecomposable subspace,

- a subspace with an unconditional basis such that no disjointly
supported block-subspaces are isomorphic,

- a subspace with an unconditional basis which is quasi-minimal.

\end{theo}

Using his game we prove:

\begin{theo} \label{minimality}
 Let $X$ be a separable Banach space. Then

(i) $X$ is ergodic
or contains a quasi-minimal subspace with an unconditional basis.

(ii) $X$ contains a perfect set of mutually totally incomparable subspaces or
a quasi-minimal subspace.

(iii) $X$ contains a perfect set of non isomorphic subspaces  or
 a block minimal subspace with an unconditional basis.
 \end{theo}

\pf First notice that because of the hereditary nature of the properties,
each of the subspaces above may be chosen to
be spanned by block-bases of a given basis. Rosendal proved that any hereditarily indecomposable Banach space $X$ is
ergodic,
and this can be proved using subspaces generated by subsequences of a basic
sequence
 in $X$ \cite{R}.
Following Bossard (who studied the particular case of a space defined by
Gowers, \cite{B}),
we may prove
that a space $X$ such that no disjointly
supported block-subspaces are isomorphic is ergodic (map
$\alpha \in 2^{\omega}$ to $[e_{2n+\alpha(n)}]_{n \in \N}$, where $(e_n)$
is the
unconditional basis of $X$). This takes care of (i).

 A space such that no disjointly
supported block-subspaces are isomorphic contains $2^{\omega}$ totally
incomparable
block-subspaces (take subspaces generated by subsequences of the basis
corresponding
to a perfect set of almost disjoint infinite subsets of $\omega$).
Also any hereditarily indecomposable space is quasi-minimal so (ii) follows.

Finally for the proof of (iii) we will first show that the statement we want
 to prove is  ${\fed{\Sigma^1_2}}$, which will be done
by showing that given a block-minimal space $X$, there is a further block subspace
$Y$ such that for $Z\leq Y$ we can find continuously in $Z$ an $X'\leq Z$ and
an isomorphism of $X'$ with $X$. The proof uses some ideas of Lopez-Abad
\cite{JLA} of coding with asymptotic sets.

By standard metamathematical facts it will then be sufficient to show the statement under Martin`s axiom
 and the negation of the continuum hypothesis. This was almost done by Bagaria and Lopez-Abad who showed it to be
consistent relative to the existence of a weakly compact cardinal, see \cite{BL}, but we will see that it can be
done in a simple manner directly from $MA+\neg CH$.

Note first that having a perfect set of non isomorphic subspaces or
containing a copy of $c_0$ are both  ${\fed{\Sigma^1_2}}$, and that if $X$
contains a copy of $c_0$ then it has a block sequence equivalent to the
unit vector basis of $c_0$, in which case the theorem holds. If on the
contrary it does not contain $c_0$ then by the solution to the distortion
problem by Odell and Schlumprecht it must contain two closed, positively
separated, asymptotic subsets of the unit sphere $A_0$ and $A_1$,
\cite{OS}. Suppose that $Y=[y]\leq X$ is block minimal. Fix a bijection $\pi$
between $\omega$ and $\Q_*^{<\om}$, the set of finite sequences of
rationals not identically zero, then for any $\alpha=
0^{n_0}10^{n_1}10^{n_2}1\ldots$ in $\ca$ there is associated a  unique
sequence $(\pi(n_0),\pi(n_1),\pi(n_2),\ldots)$ in  $(\Q_*^{<\om})^\om$.
Furthermore any element of $(\Q_*^{<\om})^\om$ gives a
unique sequence of block-vectors of
$Y$ simply by taking the corresponding finite linear combinations.

Let $D:=\mgdv (z_n)\in bb_N(X)\del  (z_n)\leq y\og \a n\; z_n\in A_0\cup A_1\og \exists^\infty n
\; z_n \in A_1\mgdh$ which is Borel in $bb_N(X)$. Then if $(z_n)\in D$ it
codes a unique infinite sequence of block vectors (not necessarily
consecutive) of $Y$, by first letting $(z_n)\mapsto \alpha\in \ca$ where
$\alpha(n)=1\hviss z_n\in A_1$ and then composing with the other coding.
Notice that this coding is continuous from $D$ to  $(\Q_*^{<\om})^\om$,
when  $\Q_*^{<\om}$ is taken discrete.

Let $E$ be the set of $(z_n)\leq y$ such that  $(z_{2n+1})\in D$ and the
function sending $(z_{2n})$ to the sequence of block-vectors of $Y$ coded
by $(z_{2n+1})$ is an isomorphism of $[z_{2n}]_{n \in \N}$ with $Y$.

$E$ is clearly Borel in $bb_N(X)$ and we claim that any block sequence
contains a further block sequence in $E$. For suppose that $z\leq y$ is
given, then we first construct a further block sequence $(z_n)$ such that
$z_{3n+1}\in A_0$ and $z_{3n+2}\in A_1$. By block minimality of $Y$ there
are a block sequence $(x_n)$ of $(z_{3n})$ isomorphic to $Y$ and a sequence
$\alpha\in \ca$ coding a sequence of block vectors $(y_n)$ of $Y$ such that
$x_n\mapsto y_n$ is an isomorphism of $[x_n]$ with $Y$ (a standard
pertubation argument shows that we can always take our $y_n$ to be a finite
rational combination on $Y$).

Now in between $x_n$ and $x_{n+1}$ there are $z_{3m+1}$ and $z_{3m+2}$, so
we can code $\alpha$ by a corresponding subsequence $(z'_n)$ of these such
that $x_n<z'_n<x_{n+1}$. The combined sequence is then in $E$. So by
Gowers' theorem  there is for any $\Delta>0$ a winning strategy $\tau$ for
II for producing blocks in $E_\Delta$ in some $Y'\leq Y$. By choosing
$\Delta$ small enough and modifying $\tau$ a bit we can suppose that the
vectors of odd index played by II are in $A_0\cup A_1$. So if $\Delta$ is
chosen small enough, a pertubation argument shows that $\tau$ is in fact a
strategy for playing blocks in $E$. By changing the strategy again we can suppose
 that II responds to block-bases in $bb_d(Y')$  by
block-bases in $bb_d(Y')$. So finally we
see that
$X$ has a block-minimal subspace iff there are $Y'=[y']$ and $Y=[y]$
with
$Y'\leq Y\leq X$ and
a continuous function  
$(f_1,f_2)=f:bb_d(Y')\til bb_d(Y')\times(\Q_*^{<\om})^\om$ such that
for all 
$z \leq y'$,
$f_2(z)$ codes a sequence $(w_n)$ of blocks of $Y$ such that
$[w_n]_{n \in \N}=Y$, and 
$f_1(z)=(v_n)\leq
z$ with $v_n\mapsto w_n$ being an isomorphism between $[v_n]_{n \in \N}$
and
$Y$.
 
The statement is therefore $\mathbf  \Sigma^1_2$, and to finish the proof we now need  the following lemma:

\begin{lemm}($MA_{\sigma-centered}$) Let $A\subset bb_\Q$ be linearly ordered under $\leq^*$ of
cardinality strictly less than the continuum, then there is an
$y_{\infty}\in bb_\Q$ such that $x_0\leq^* y$ for all $y\in A$.
\end{lemm}

\pf For $s \in fbb_{\Q}$ and $y \in bb_{\Q}$, denote by
$(s,y)$ the set of block-bases in $bb_{\Q}$ of the form $s^{\frown}z$ for $z \leq y$.
Let $\poset=\mgdv
(s,y)\del s\in fbb_\Q \og y\in A \mgdh$, ordered by the inclusion.  As a preliminary remark,
 note
that if $(t,z) \subset (s,y)$, then $s\preceq t$, $t \setminus s \leq y$, and $z \leq^* y$;
conversely, if $s \in fbb_{\Q}$ and $z \leq^* y$, then extensions $t$ of $s$,
with $t\setminus s \leq z$ far enough, are such that $(t,z) \subset (s,y)$. 

Put $D_n=\mgdv  (s,y)\in \poset\del |s|\geq n\mgdh$ and  $D_y=\mgdv  (t,z)\in
\poset \del z \leq^* y\mgdh$. Then $D_n$ and $D_y$, for $y\in A$, are dense in $\poset$, i.e.,
any element in $\poset$ has a minorant in $D_n$ (resp. $D_y$). To see that $D_n$ is dense, just take for any
given $(s,y)\in \poset$ some extension $s'$ of $s$ such that $s'\setminus s\leq y$ and
$|s'|\geq n$, then  $(s',y)\in D_n$ and $(s',y)\subset (s,y)$. On the
other hand to see that $D_y$ is dense for $y\in A$ suppose $ (s,z)\in \poset$ is given, then
as $A$ is linearly ordered by $\leq^*$, let $w$ be the minimum of $z$ and $y$. By the preliminary
remark, $(s',w) \subset (s,z)$ for a long enough extension $s'$ of $s$ such that $s' \setminus s \leq w$,
and
 as $w
\leq^* y$,
$(s',w)$ is in $D_y$.

Let $\poset_s=\mgdv (s,y)\del y\in A\mgdh$, which is centered in $\poset$, i.e.,
every finite subset of $\poset_s$ has a common minorant in $\poset$. This follows
 from the same argument as above, using the preliminary remark.
So since $s$ is supposed to be rational, we see that $\poset$ is $\sigma$-centered, i.e., a countable
union of centered subsets. Notice that as $|A|<2^\om$, there are less than continuum many dense sets
$D_n$ and $D_y$. So by $MA_{\sigma-centered}$ there is a filter $G$ on $\poset$ intersecting each of
these sets.

Suppose that $(s,y)$ and $(t,z)\in G$ then as $G$ is a filter, they have a common minorant
$(v,w)\in G$, but then $s\preceq v$ and $t\preceq v$, so either $s\preceq t$ or $t\preceq s$.
Therefore $y_{\infty}:=\for\mgdv s\in fbb_\Q\del \exists y\; (s,y)\in
G\mgdh$ is a block-basis. Furthermore as $G$ intersects all of $D_n$ for
$n<\om$ we see that $y_{\infty}$ is an infinite block-basis.

We now prove that $y_{\infty} \leq^* y$ for all $y \in A$. Since $G$
intersects
$D_y$, without loss of generality
 we may
assume
 that
$(s,y)\in G$ for some $s$. Then
$y_{\infty}\setminus s\leq y$. For if $t\preceq y_{\infty}$ and $(t,z)\in
G$, take
$(v,w)\in G$ such that
 $(v,w)\subset (t,z)$ and  $(v,z)\subset (s,y)$,
 then $s,t\preceq v\preceq y_{\infty}$ and $t\setminus s\preceq v\setminus
s\leq y$, and therefore as $t$ was arbitrary $y_{\infty}\setminus s\leq
y$.\pff  

Suppose now that $X$ does not have a perfect set of non isomorphic subspaces,
 then by Burgess' theorem (see \cite {kec} (35.21)), it has at most $\om_1$ isomorphism classes of
subspaces, and in particular as we are supposing the continuum hypothesis not to hold, less than continuum many.
Let
$(X_\xi)_{\xi<\om_1}$ be an enumeration of an element from each class. Then if none of these are minimal, we can
construct inductively a $\leq^*$ decreasing sequence $(Y_\xi)_{\xi<\om_1}$ of rational block-subspace such
that
$X_\xi$ does not embed into $Y_\xi$ and using the above lemma find some $Y_{\om_1}$ diagonalizing the whole
sequence. By taking, e.g., the subsequence consisting of every second term
of the basis of
$Y_{\om_1}$ one can suppose that 
$Y_{\om_1}$ embeds into every term of the sequence   $(Y_\xi)_{\xi<\om_1}$ and that therefore in particular  
$Y_{\om_1}$ is isomorphic to no   $X_\xi$, $\xi<\om_1$, which is impossible. This finishes the proof of
the theorem.
\pff

We remark that in the case of $X$ not containing a minimal subspace. then there is in fact a
 perfect set of subspaces such each two of them do not both embed into each other,
 which is slightly stronger than
saying that they are non isomorphic.

\

  We
recall our result from
\cite{FR},\cite{R} in a slightly modified form.

\begin{theo} \label{theocv}
Let $X$ be a Banach space with a basis $\{e_i\}$. Then $E_0$
Borel reduces to isomorphism on
subspaces spanned by subsequences of the basis, or there exists
 a sequence $(F_n)_{n \geq 1}$ of successive finite
 subsets of $\N$ such that for any infinite subset $N$ of $\N$,
 if  $N \cap [\min(F_n),\max(F_n)]=F_n$ for infinitely many $n$'s,
  then
 the space $[e_i]_{i \in N}$ is isomorphic to $X$. It follows that if $X$
is non ergodic with an unconditional basis, then it
is isomorphic to its hyperplanes, to its square, and more generally to $X
\oplus Y$ for any subspace $Y$ spanned by
a subsequence of the basis.
 \end{theo}

Indeed, by \cite{FR}, improved in \cite{R}, either $E_0$ Borel reduces to
isomorphism on
subspaces spanned by subsequences of the basis, or the set of infinite
subsets of $\N$ spanning a space
isomorphic to $X$ is residual  in $2^{\omega}$; the characterization in
terms of finite subsets of $\N$
 is then a
classical characterization of residual  subsets of $2^{\omega}$ (see \cite{kec},
 or the remark at the end of Lemma 7 in \cite{FR}). Both proofs are similar to (and simpler than)
the following proofs of Proposition \ref{prop5} and Proposition \ref{prop6} for block-subspaces.
 The last part of
the theorem is specific to the case of subspaces spanned by subsequences and is also proved
 in \cite{FR}.

\

We now wish to extend this result to the set of block-bases, for which it is useful to use the
 Polish space $bb_d(X)$. Unless stated otherwise this is the topology refered to. 

As before, the notation $x=(x_n)_{n \in \N}$ will be used to denote
an infinite block-sequence; $\tilde{x}$ will denote a
finite block-sequence, and $|\tilde{x}|$ its
length as a sequence, $supp(\tilde{x})$ the union of the supports of the
terms of $\tilde{x}$. For two finite block-sequences $\tilde{x}$ and
$\tilde{y}$, write $\tilde{x}<\tilde{y}$ to mean that they are
successive. For a sequence of successive finite block-sequences
$(\tilde{x}_i)_{i
\in I}$, we denote the concatenation of the block-sequences by
$\tilde{x}_1^{\frown} \ldots^{\frown} \tilde{x}_n$ if the sequence is finite or
$\tilde{x}_1^{\frown} \tilde{x}_2^{\frown}\ldots$ if it is infinite, and
we denote by $supp(\tilde{x}_i, i \in I)$ the support of the
concatenation, by
$[\tilde{x}_i]_{i
\in I}$ the closed linear span of the concatenation.  For a finite block
sequence
$\tilde{x}=(x_1,\ldots,x_n)$, we denote by
$N(\tilde{x})$
the set of elements of $bb_d(X)$ whose first $n$ vectors are
$(x_1,\ldots,x_n)$.

\begin{prop} \label{prop5}
Let $X$ be a Banach space with a Schauder basis. Then either $X$ is
ergodic, or there exists
$K \geq 1$ such that a residual  set of bloc-sequences in
$bb_{d}(X)$ span spaces
mutually $K$-isomorphic.
\end{prop}

\pf
  The
relation of isomorphism is either meagre or non-meagre in $bb_d(X)^2$.
 First assume that it
is meagre.
Let $(U_n)_{n \in \N}$ be a decreasing sequence of dense open subsets of
$bb_d(X)^2$
so that $\cap_{n \in \N} U_n$ does not intersect $\simeq$.
We build by induction successive finite blocks $\{\tilde{a}_n^0, n \in \N\}$
 and $\{\tilde{a}_n^1, n \in \N\}$
such that for all $n$, $|\tilde{a}_n^0|=|\tilde{a}_n^1|$, and
$supp(\tilde{a}_n^i) < supp(\tilde{a}_{n+1}^j)$ for all $(i,j) \in 2^2$.
For $\alpha \in 2^{\omega}$, we let $x(\alpha)$ be the concatenated
infinite block-sequence
${\tilde{a}_0^\alpha(0)}^{\frown} {\tilde{a}_1^\alpha(1)}^{\frown}\ldots$;
 for $n \in \N$ and $\beta \in 2^{n}$,
we let $\tilde{x}(\beta)$ be the concatenated finite block-sequence
${\tilde{a}_0^\alpha(0)}^{\frown} \ldots^{\frown}
\tilde{a}_{n-1}^{\alpha(n-1)}$.
We require furthermore from the
sequences $\{\tilde{a}_n^0\}$ and $\{\tilde{a}_n^1\}$
that for each $n \in \N$, each $\beta$ and $\beta'$ in $2^n$,
$$N(\tilde{x}(\beta^{\frown}0)) \times N(\tilde{x}({\beta'}^{\frown}1))
\subset U_n.$$
Before explaining the construction, let us check that with these conditions,
the map $\alpha \rightarrow x(\alpha)$ realizes a Borel reducing of $E_0$
into
$(bb_d(X),\simeq)$. Indeed, when $\alpha E_0 \alpha'$, the corresponding
 sequences
differ by at most finitely many vectors, and since we took care that
$|\tilde{a}_n^0|=|\tilde{a}_n^1|$ for all $n$, $x(\alpha)$ and
$x(\alpha')$ span isomorphic subspaces. On the other hand, when
$\alpha$ and  $\alpha'$ are not $E_0$-related, without loss of generality
  there is an infinite set $I$ such
that for all $i \in I$, $\alpha(i)=0$ and $\alpha'(i)=1$; it follows that
 for all $i \in I$,
$(x(\alpha),x(\alpha'))$ belongs to $U_i$, and so by choice of the $U_n's$,
$(x(\alpha),x(\alpha'))$ does not belong to $\simeq$.

Now let us see at step $n$ how to construct the sequences:
given a pair $\beta_0,\beta'_0$ in $(2^n)^2$, using the fact that $U_n$ is
dense and open, the pair
$\tilde{x}(\beta_0),\tilde{x}(\beta'_0)$  may be extended to a pair of finite
 successive block-sequences of the form
$(\tilde{x}(\beta_0)^{\frown}
\tilde{z}_0,\tilde{x}(\beta'_0)^{\frown} \tilde{z}'_0)$ such that
$N(\tilde{x}(\beta_0)^{\frown}
\tilde{z}_0)\times N(\tilde{x}(\beta'_0)^{\frown} \tilde{z}'_0)) \subset
U_n$, and we may require that
$supp(\tilde{x}(\beta_0)) \cup supp(\tilde{x}(\beta'_0)) <
supp(\tilde{z_0}) \cup supp(\tilde{z}'_0)$. Given an other pair
 $\beta_1,\beta'_1$ in $(2^n)^2$, the pair
 $(\tilde{x}(\beta_1)^{\frown} \tilde{z_0},\tilde{x}(\beta'_1)^{\frown}
\tilde{z}'_0)$
may be extended to a pair of finite successive block-sequences
$(\tilde{x}(\beta_1)^{\frown}
\tilde{z}_1,\tilde{x}(\beta'_1)^{\frown} \tilde{z}'_1)$ such that
$N(\tilde{x}(\beta_1)^{\frown}
\tilde{z}_1)\times N(\tilde{x}(\beta'_1)^{\frown} \tilde{z}'_1) \subset
U_n$. Here with our notation
$\tilde{z}_1$ extends $\tilde{z}_0$ and
$\tilde{z}'_1$ extends $\tilde{z}'_0$.  Repeat this $(2^n)^2$ times
to get $\tilde{z}_{4^n-1},\tilde{z}'_{4^n-1}$ such that
for all $\beta$ and $\beta'$ in $2^n$,
$$N(\tilde{x}(\beta)^{\frown}\tilde{z}_{4^n-1}) \times
N(\tilde{x}(\beta')^{\frown}
\tilde{z}'_{4^n-1}) \subset U_n.$$
Finally extend $(\tilde{z}_{4^n-1},\tilde{z}'_{4^n-1})$ to
$(\tilde{a}_n^0,\tilde{a}_n^1)$ such that $|\tilde{a}_n^0|=|\tilde{a}_n^1|$; we
still have for all $\beta$ and $\beta'$ in $2^n$,
$N(\tilde{x}(\beta)^{\frown}\tilde{a}_n^0) \times
N(\tilde{x}(\beta')^{\frown}\tilde{a}_n^1) \subset U_n,$ i.e.
with our notation,
$$N(\tilde{x}(\beta^{\frown}0)) \times N(\tilde{x}({\beta'}^{\frown}1))
\subset U_n.$$

\

 Now assume the relation of isomorphism is non-meagre in $bb_d(X)^2$.
As the relation is analytic it has the Baire
property and $bb_d(X)$ is Polish, so by Kuratowski-Ulam (Theorem
8.41 in \cite{kec}), there
must be some non-meagre section, that is, some isomorphism class $\cal A$ is
non-meagre. Fix a block-sequence
$x$ in this class, then clearly, for some constant
$C$, the set ${\cal A}_C$ of blocks-sequences spanning a space
$C$-isomorphic to
$[x]$ is non-meagre.
Now being analytic, this set has the Baire property, so is residual in some
basic
open set $U$, of the form $N(\tilde{x})$, and we may assume $n \geq 1$.

We now prove that ${\cal A}_k$ is residual  in $bb_d(X)$ for
$k=Cc(2\max(supp(\tilde{x})))$; The conclusion of the proposition then holds
for
$K=k^2$.
Recall that for $n \in \N$, $c(n)$ denotes a constant such for any Banach space
$X$, any $n$-codimensional subspaces of $X$ are $c(n)$-isomorphic (Lemma 3
in \cite{FR}).
So
let us assume
 $V=N(\tilde{y})$ is some basic open set in $bb_d(X)$ such that
${\cal A}_k$ is meagre in $V$. We may assume that $|\tilde{y}|>|\tilde{x}|$
and write $\tilde{y}=\tilde{x}'^{\frown}\tilde{z}$ with
$\tilde{x}<\tilde{z}$ and $|\tilde{x}'| \leq \max(supp(\tilde{x}))$.
Choose $\tilde{u}$ and $\tilde{v}$ to be finite sequences of blocks such
that $\tilde{u},\tilde{v}>\tilde{z}$, $|\tilde{u}|=|\tilde{x}'|$ and
$|\tilde{v}|=|\tilde{x}|$, and such that 
$\max(supp(\tilde{u}))=max(supp(\tilde{v}))$.
Let $U'$ be the basic open set
$N(\tilde{x}^{\frown}\tilde{z}^{\frown}\tilde{u})$ and
let $V'$ be the basic open set
$N(\tilde{x}'^{\frown}\tilde{z}^{\frown}\tilde{v})$.
Again ${\cal A}_C$ is residual  in $U'$ while
${\cal A}_k$ is meagre in $V'$.

Now let $T$ be the canonical map from $U'$ to $V'$.
For all $u$ in $U'$, $T(u)$ differs from at most
 $|\tilde{x}|+\max(supp(\tilde{x})) \leq
2\max(supp(\tilde{x}))$ vectors from $u$, so $[T(u)]$ is
$c(2\max(supp(\tilde{x})))$ isomorphic to $[u]$. Since
$k=Cc(2\max(supp(\tilde{x})))$ it
follows that
${\cal A}_{k}$ is residual  in $V' \subset V$ . The contradiction
follows by
choice of
$V$. \pff

\

By analogy with the definition of atomic measures,
we may see Proposition
\ref{prop5} as stating that a non-ergodic Banach spaces with a basis
must be "atomic" for its block-subspaces.

We now want to give a characterization of residual  subsets of $bb_d(X)$.
If $\cal A$ is a subset of $bb_d(X)$ and $\Delta=(\delta_n)_{n \in \N}$ is a
sequence of
 positive real numbers, we denote by ${\cal A}_{\Delta}$ the usual
$\Delta$-expansion of
$\cal A$ in $bb_d(X)$, that is $x=(x_n) \in {\cal A}_{\Delta}$ iff
there exists $y=(y_n) \in {\cal A}$ such that $\norm{y_n-x_n} \leq
\delta_n, \forall n \in
\N$.
Given a finite block-sequence $\tilde{x}=(x_1,\ldots,x_n)$, we say that a
(finite or infinite)
block-sequence $(y_i)$ {\em passes through} $\tilde{x}$ if
there exists some integer $m$ such that $\forall 1 \leq i \leq n$,
$y_{m+i}=x_i$.

\begin{prop} \label{prop6}
Let $\cal A$ be residual  in $bb_d(X)$. Then
 for all $\Delta>0$, there exist successive finite block-sequences
  $(\tilde{x}_n), n \in \N$ such that any element of $bb_d(X)$ passing trough
  infinitely many of the $\tilde{x}_n$'s is in ${\cal A}_{\Delta}$.
  \end{prop}

\pf
Let $(U_n)_{n \in \N}$ be a sequence of dense open sets, which we may
assume to be decreasing,
such that $\cap_{n \in \N} U_n \subset {\cal A}$. Without loss of
generality we may also
assume $\Delta$ to be decreasing. In the following, block-vectors are
always taken in $\mathbf {Q}$,
in the intention of building elements of $bb_d(X)$.

First $U_0$ is open so there exists $\tilde{x}_0$ a fi\-ni\-te block
se\-quen\-ce
such that $N(\tilde{x}_0) \subset U_0$.
Now let us choose some $n_1>max(supp(\tilde{x}_0))$ and let us take an
arbi\-trary
block-vector $z_1$, and let $N_1=\min(supp(z_1))$.
Let $F_{<1}$ be a fi\-nite set of finite block sequences forming an
$\delta_{N_1}$- net
for all finite block sequences supported before $N_1$ and let $F_{01}$ be
a finite set of finite block sequences forming an $\delta_{N_1}$- net
for all finite block sequences supported after $\tilde{x}_0$ and before
$N_1$. Let $G_1=\{\tilde{x}_0^{\frown}\tilde{y}, \tilde{y} \in F_{01}\}$ and let
$F_1=F_{<1} \cup G_1$.
Using the fact that $U_1$ is dense open,
we may construct successively a finite block sequence $\tilde{x}_1$  which
prolonges
$z_1$, so that $\min(supp(\tilde{x}_1))=N_1>\max(supp(\tilde{x}_0))$,
and
such that for
any $\tilde{f}_1 \in F_1$,
$N(\tilde{f}_1^{\frown} \tilde{x}_1)$ is a subset of $U_1$.

Let us now write what happens at the $k$-th step.
We choose some $N_k>\max(supp(\tilde{x}_{k-1}))$ and  an arbitrary
block $z_k$ whose support starts at $N_k$.
We let $F_{<k}$ be a finite set of finite block sequences forming an
$\delta_{N_k}$- net
for all finite block sequences supported before $n_k$ and for all $i<k$,
 we let $F_{ik}$ be
a finite set of finite block sequences forming an $\delta_{N_k}$- net
for all finite block sequences supported after $\tilde{x}_i$ and before
$N_k$. For any $I=\{i_1<i_2<\cdots<i_m=k\}$, we let $G_{I}$ be the set of
finite block sequences $\tilde{z}$ passing through every $i$ in $I$,
such that
the finite sequence of blocks
of $\tilde{z}$ supported before $\tilde{x}_{i_0}$
is in $F_{<i_0}$ and such that for all $j<m$, the finite sequence of blocks
of $\tilde{z}$
 supported between
$\tilde{x}_{i_j}$ and $\tilde{x}_{i_{j+1}}$ is in $F_{i_j i_{j+1}}$.
And we let $F_k$ be the union of all $G_I$ over all possible subsets
of
$\{1,2,\ldots,k\}$ containing $k$.
 Using the fact that $U_k$ is dense open,
we may construct successively a finite block sequence $\tilde{x}_k$  which
prolonges
$z_k$, so that $\min(supp(\tilde{x}_k))=N_k>\max(supp(\tilde{x}_{k-1}))$,
and
such that for
any $\tilde{f}_k \in F_k$,
$N(\tilde{f}_k^{\frown} \tilde{x}_k)$ is a subset of $U_k$.

Repeat this construction by induction, and now let
 $z$ be a block-sequence passing through $\tilde{x}_n$
 for $n$ in an infinite set $\{n_k, k
\in \N\}$. We may write $z=\tilde{y}_0^{\frown} \tilde{x}_{n_0}^{\frown}
\tilde{y}_1^{\frown} \tilde{x}_{n_1}^{\frown} \ldots$, where
 $\tilde{y}_0$ is supported before $\tilde{x}_{n_0}$ (we may assume that
$n_0>0$)
  and for $k>0$,
$\tilde{y}_k$ is supported between $\tilde{x}_{n_{k-1}}$ and
$\tilde{x}_{n_k}$.

Let $\tilde{f}_0 \in F_{<n_0}$ be $\delta_{N_{n_0}}$ distant from
$\tilde{y}_0$, and
for any $k>0$, let $\tilde{f}_k \in F_{n_{k-1}n_k}$ be
$\delta_{N_{n_k}}$ distant from $\tilde{y}_k$. Then it is clear that
$z$ is $\Delta$ distant from
$f=\tilde{f}_0^{\frown}\tilde{x}_{n_0}^{\frown}\tilde{f}_1^{\frown}\tilde{x}
_{n_1}
^{\frown} \dots$.
Indeed, consider a term $z_n$ of the block-sequence $z$: if it appears as a
term of some finite sequence $\tilde{x}_{n_k}$
then its distance to the corresponding block $f_n$ of $f$ is $0$; if it appears
as a term of  some $\tilde{y}_k$ then it is
less than $\delta_{N_k}$-distant from the
block $f_n$, and $N_k>\max(supp(z_n)) \geq n$, so it is less
than $\delta_n$-distant from $f_n$.

It remains to check that $f$ is in $\cal A$.
But for all $K$, the finite sequence $\tilde{g}_K=\tilde{f}_0^{\frown}
\tilde{x}_{n_0}^{\frown}
\ldots^{\frown} \tilde{f}_k^{\frown} \tilde{x}_{n_k}$
is an element of $G_{\{n_1,\ldots,n_K\}}$ so is an element of
$F_K$; it follows that $N(\tilde{g}_K)$
is a subset of $U_{n_K}$ and so that $f$ is in $U_{n_K}$.
Finally $f$ is in $\cap_{k \in \N} U_{n_k}$ so is in $\cal A$.  \pff

\

Conversely,
given successive blocks
$\tilde{x}_n$, the set of block-sequences passing through infinitely many
 of the $\tilde{x}_n$'s is residual:
  for a given $\tilde{x}_n$,
"$(y_k)_{k \in \N}$ passes through $\tilde{x}_n$" is open and
 "$(y_k)_{k \in \N}$ passes through infinitely many of the $\tilde{x}_n$'s"
is equivalent to
 "$\forall m \in \N, \exists n>m \in \N$: $(y_k)$ passes through
$\tilde{x}_n$", so is
$G_{\delta}$, and clearly dense. If the set $\cal A$ considered is an
isomorphism class,
 then it is invariant under small
enough $\Delta$-perturbations, and  so we get an equivalence:
$\cal A$ is residual iff there exist successive finite block-sequences
  $(\tilde{x}_n), n \in \N$ such that any element of $bb_d(X)$ passing trough
  infinitely many of the $\tilde{x}_n$'s is in ${\cal A}$.

  Finally, as any element of $bb(X)$ is arbitrarily close to an element of $bb_d(X)$,
 the following theorem holds:

\begin{theo} \label{theo7}
 Let $X$ be a Banach space with a basis. Then either $X$ is ergodic
 or there exists  $K \geq 1$, and a sequence
of successive finite block-sequences $\{\tilde{x}_n\}$ such that all
 block-sequences passing through infinitely many of
the $\{\tilde{x}_n\}$'s span mutually $K$-isomorphic subspaces.
\end{theo}

\

If in addition the basis is unconditional, then we may use the projections
to get
 further properties
of the residual class.

 \begin{coro} \label{coro8} Let $X$ be a non ergodic Banach space with an
 unconditional basis.  Denote by $A$ an element of the residual class
 of isomorphism in $bb_d(X)$.
  Then for any block-subspace $Y$ of $X$, $A \simeq A \oplus Y$.
If
$X$ is hereditarily block-minimal, then all residual classes in $bb_d(Y)$,
for block-subspaces
$Y$ of $X$,
are isomorphic.
\end{coro}

 \pf Let $\{e_i\}$ be the unconditional basis of $X$ and let
$\{\tilde{x}_n\}$ be given by
  Theorem \ref{theo7}. Consider an arbitrary block-subspace $Y$ of $X$.
  Its natural basis is unconditional and $Y=[y_i]_{i \in \N}$ is not ergodic as well.
Let, by Theorem \ref{theocv},
$(F_n)_{n \geq 1}$ be successive finite
 subsets of $\N$ such that for any infinite subset $N$ of $\N$,
 if  $N \cap [\min(F_n),\max(F_n)]=F_n$ for infinitely many $n$'s,
  then
 the space $[y_i]_{i \in N}$ is isomorphic to $Y$.
 Passing to subsequences we may assume that
 for all $n$ in $\N$,
 $\tilde{x}_n < \cup_{i \in F_n} supp(y_i)< \tilde{x}_{n+1}$.
 Then
 $$A \simeq [\tilde{x}_n]_{n \in \N} \oplus [y_i]_{i \in \cup_{n \in \N}
F_n} \simeq A \oplus Y.$$
 If now $X$ is hereditarily block-minimal, and $B$ belongs to the
resi\-dual class in $bb_d(Y)$, for $Y \leq X$,
 then by the abo\-ve $A \simeq A \oplus B$; but also $A$ is block-minimal
so some
copy of $A$ embeds as a block-subspace of $Y$, so $B \simeq B \oplus A$. \pff

\

 Consider the property that every block-subspace $Y$ satisfies
 $A \simeq A \oplus Y$.
  We may think
 of this property as an algebraic property characterizing large subspaces
 in the sense that a large subspace
 should intuitively "contain" other subspaces, and more
  importantly, a space
 should have at most one large subspace (here if $A$ and $A'$ satisfy the
property,
 $A \simeq A \oplus A' \simeq A'$!).
  Notice that as $X$ is not ergodic,
 all block-subspaces are isomorphic to their squares by \cite{FR}
  and so the property above
  is equivalent to saying
  that every block of $X$ embeds complementably in $A$ (i.e.
  $A$ is {\em complementably universal for
  $bb(X)$}). Generally, a space $A$ is said to be  complementably universal
for a
  class $\cal C$ of Banach spaces if every element of $\cal C$
  is isomorphic to a complemented subspace of $X$.
  It is known that no separable Banach space is complementably universal for the
  class of all separable Banach spaces (\cite{LT}, Th. 2.d.9), but there
exists a Banach space $X_U$
  with an unconditional basis which is complementably universal for the class of
  all Banach spaces with an unconditional basis (\cite{P}), and so for
$bb(X_U)$ in
  particular.

  Combining Theorem \ref{minimality} and Corollary \ref{coro8}, we get

  \begin{theo} \label{theo7.5}
 Any Banach space is ergodic or contains a subspace with an unconditional
basis which is
 complementably universal for the family of its block-subspaces.
 \end{theo}

 We now study this property in more details. We also
  see how Theorem \ref{theo7}
 may be used to obtain uniformity results.

\begin{defi} \label{defi9}
 Let $X$ and $Y$ be Banach spaces such that $X$ has a Schauder basis.
  $Y$ is said to be complementably universal for $bb(X)$
if every block subspace of $X$ embeds complementably in $Y$.
\end{defi}

\begin{lemm} \label{decomposability}
Let $X$ be a Banach space. Any Banach space complementably universal for
$bb(X)$ is decomposable.
\end{lemm}

\pf  Let $A$ be complementably universal for $bb(X)$ and indecomposable.
First note
that $X$ embeds complemen\-tably in $A$, so must be isomorphic to
 a finite-codimensional
subspace of $A$. As well, any
block subspace of $X$
is isomorphic to a finite-codimen\-sional subspace of $A$ and  so none of
them is
decomposable either. It follows easily that no subspace of $X$ is decomposable.
In other words, $X$ is heredi\-tarily indecomposable. It follows also
that $X$ is isomorphic to a proper (infinite-dimensional) subspace, and
this is a contradiction with properties of hereditarily indecomposable
spaces. \pff

\

 To quantify the property of complementable universality, let us
 define
$dec_X(Y)=\inf KK'$, where the infimum runs over all
 couples $(K,K')$
such that $Y$ is $K$-isomorphic to a $K'$-complemented subspace of $X$.
Of course, $dec_X(Y)=+\infty$ iff $Y$ does not embed comple\-mentably in
 $X$. We shall say that a space $A$ is {\em $C$-complementably universal}
for $bb(X)$
  if every block-sub\-space of $X$
is
$K$-isomorphic to some $K'$-complemented sub\-space of $A$, for some $K$
and $K'$
such that $KK' \leq C$, that is, if $\sup_{Y \leq X} dec_A(Y) \leq C$.

\begin{lemm} \label{lemm11} Assume $A$, $B$, $C$ are Banach spaces with
bases. Then
$$dec_A(C) \leq dec_A(B)^2 dec_B(C).$$
\end{lemm}

\pf Let $\e$ be positive. Let $P_{B}$ be a projection defined on $B$ and
 $\alpha_{BC}$ be an isomorphism from $P_B(B)$ onto $C$ such that
  $\norm{P_{B}}.\norm{\alpha_{BC}}.\norm{\alpha_{BC}^{-1}} \leq dec_B(C)+\e$.
 Let $P_{A}$ be a projection defined on $A$ and $\alpha_{AB}$ be an isomorphism
 from $P_A(A)$ onto $B$ such that
$\norm{P_{A}}.\norm{\alpha_{AC}}.\norm{\alpha_{AB}^{-1}} \leq dec_A(B)+\e$.

 We  let $P=\alpha_{AB}^{-1} P_B \alpha_{AB} P_A$, defined on $A$; it is
easily checked that $P$
 is a projection. We let $\alpha=\alpha_{BC}\alpha_{AB}$: it is an
isomorphism from $P(A)$ onto $C$.
 Then $$dec_A(C) \leq \norm{P}.\norm{\alpha}.\norm{\alpha^{-1}} \leq
 (dec_A(B)+\e)^2 (dec_B(C)+\e). $$ \pff

 \begin{prop} \label{prop12} Let $X$ be a Banach space with
 a Schauder basis and let $A$ be complementably universal for $bb(X)$.
 Then there exists $C \geq 1$ such that every finite dimensional
  block-subspace of $X$ $C$-embeds
 complementably in $A$.
 \end{prop}

 \pf  First it is clear that it is enough to restrict ourselves to elements
of $bb_d(X)$ with the previously defined topology.
We let for $k \in \N$, ${\cal A}_k$ denote the
set of block-subspaces of $X$ which are $k$-isomorphic to
some $k$-complemented subspace of $A$. Now it
is clear that one of the ${\cal A}_k$ must
be non-meagre. This set is analytic, so has the Baire property,
 so is residual in some
basic
open set $U$, of the form $N(\tilde{u})$. 
We now show that ${\cal A}_{K}$ is residual for
$K=kc(2\max(supp(\tilde{u}))$. Otherwise, as in Proposition \ref{prop5},
we may assume
${\cal A}_K$ is meagre in $V=N(\tilde{y})$, and
${\cal A}_k$ is residual
in $U'=N(\tilde{x})$ where $\tilde{x}$ extends $\tilde{u}$,
 $|\tilde{x}| \leq
2\max(supp(\tilde{u}))$ and $\max(supp(\tilde{x}))=\max(supp(\tilde{y}))$.

Now let $T$ be the canonical map from $U'$ to $V$.
For all $u$ in $U'$, $T(u)$ differs from at most $q \leq
 2\max(supp(\tilde{u}))$ vectors from $u$,
so the space $[T(u)]$ is $c(2\max(supp(\tilde{u}))$ isomorphic to $[u]$.
So $T(u)$ is in ${\cal A}_K$ whenever $u$ is in ${\cal A}_k$.
It
follows that
${\cal A}_{K}$ is residual in $V$, a contradiction.

Now consider any finite dimensional space $F$ generated by
a finite block-sequence of $A$, $F=[x_1,\ldots,x_p]$:
 it may be extended
to a block-sequence $x=(x_i)_{i \in \N}$ in ${\cal A}_K$, that is
$dec_A([x]) \leq K^2$. But also $dec_{[x]}(F) \leq c$, where
$c$ is the constant of the basis, so by Lemma \ref{lemm11},
$dec_A(F) \leq cK^4$. \pff

\begin{prop} \label{prop13} Let $X$ be a space with an unconditional basis.
If $A$ is complementably universal for $bb(X)$ and isomorphic to its
square, then
$A$ is $C$-complementably universal for $bb(X)$ for some $C \geq 1$.
\end{prop}

\pf The first part of the proof is as above to get $K \in \N$ such that
${\cal A}_K$, the set of block-subspaces of $X$ which are $K$-isomorphic to
some $K$-complemented subspace of $A$, is residual.
So by Proposition \ref{prop6}, there exists a sequence ${\tilde{x}_n}$ of
successive finite blocks
  such that any
block passing through infinitely many of the $\tilde{x}_n's$ is in ${\cal
A}_{2K}$.

Let now $Y=[y_n]_{n \in \N}=[y]$ be an arbitrary block-subspace of $X$.
We may define a sequence $(\tilde{y}_i)$ of finite block-sequences
with $\tilde{y}_1^\frown \tilde{y}_2^\frown \ldots=y$
 and
a subsequence of $\{\tilde{x}_n\}$, denoted $\{\tilde{x}_i\}$, such that for
all
$i$,
$$supp(\tilde{y}_{i-1})<supp(\tilde{x}_{i})<supp(\tilde{y}_{i+1}).$$

We let $w=\tilde{y}_1^\frown \tilde{x}_{2}^\frown \tilde{y}_3^\frown \ldots$ and
 $w'=\tilde{x}_1^\frown \tilde{y}_{2}^\frown \tilde{x}_3^\frown \ldots$
It is clear that $Y_1=[\tilde{y}_{2i-1}]_{i \in \N}$ is $c$-comple\-mented
in $[w]$,
 where $c$ is the con\-stant of un\-condi\-tiona\-lity of the basis;
 so  $dec_{[w]}(Y_1) \leq c$. But we know that
$dec_A ([w]) \leq 4 K^2$, so by Lemma \ref{lemm11}, $dec_A(Y_1) \leq 16 K^4 c$.
Likewise if we denote $[\tilde{y}_{2i}]_{i \in \N}$ by $Y_2$ and use $[w']$,
we prove that $dec_A (Y_2) \leq 16 K^4 c$.
It follows that $Y_1 \oplus Y_2$ satisfies
$dec_{A \oplus_1 A} (Y_1 \oplus Y_2) \leq (16 K^4 c)^2$ (here $\oplus_1$
denotes the
$l_1$-sum) , and, if $D$ is such that
$A$ is $D$ isomorphic to its square $A \oplus_1 A$, that $dec_A(Y) \leq 2^9
D^2 K^8 c^3$.
\pff

In view of the fact that by Theorem \ref{theocv}, any unconditional basic
sequence
in a non-ergodic Banach space spans a
space isomorphic to its square, the previous proposition may be
applied to Theorem
\ref{theo7.5}: a non-ergodic Banach space contains a subspace $X$ with an
unconditional basis which is
{\em uniformly} complementably universal  for $bb(X)$.

\

Consider now spaces with a basis with the stronger property that every
block-subspace
is comple\-mented.
It is well-known that every block-subspace of $l_p$ or $c_0$ is complemented,
and the same is true for Tsirelson's spaces $T_{(p)}$ (see \cite{CS}) or
for spaces $(\sum_{n=1}^{+\infty} \oplus l_s^n)_p$, the relevant case being
$s \neq p$ (see \cite{LT}). All these
 examples
are asymptotically $l_p$ or $c_0$, and using results of Maurey, Milman and
Tomczak-Jaegermann (\cite{MMT}), we shall prove that this is not by chance.

\

We recall the definition of an asymptotically $l_p$ space with a basis.
Consider the so-called {\em asymptotic game} in $X$, where
 Player I plays integers $(n_k)$ and Player II
plays successive unit vectors $(x_k)$ in $X$ such that $supp(x_k)>n_k$ for
all
$k$. Then $X$ is asymptotically $l_p$ if there exists a constant $C$ such
that for any $n \in \N$, Player I has a winning strategy in the asymptotic
game of length $n$ for producing a sequence $C$-equivalent to the unit basis
of $l_p^n$.

\begin{prop} \label{uniformity of complementation} Let $X$ be a space with
an unconditional basis
 such that
every block-subspace of $X$ is complemented. Then there exists $C>0$ such that
every block-subspace of $X$ is $C$-complemented.
\end{prop}

\pf  Once again, we may and  do restrict ourselves to block-bases in $bb_d(X)$ with our usual topology.
Let $c$
denote the unconditional constant of the basis of
$X$. We note a fact: if an infinite block-sequence $x=(x_n)_{n \in \N}$
 (resp. $x'=(x'_n)_{n \in \N}$) spans a
$C$-complemented (resp. $C'$-complemented)
subspace of $X$, and if $x$ and $x'$ are disjointly supported, then
$[x] \oplus [x']$ is $c^2(C+C')$-complemented in $X$.
Indeed, if $q$ (resp. $q'$) is the canonical projection onto
$\cup_{n} supp(x_n)$ (resp. $\cup_{n} supp(x'_n)$), and if $p$ (resp. $p'$) is
a projection onto $[x]$ of norm less than $C$ (resp. onto $[x']$ of
norm less than $C'$),
then $qpq+q'p'q'$ is a projection onto $[x] \oplus [x']$ of norm less than
 $c^2(C+C')$.

It follows that if $[x]$ is $C$-complemented in $X$, then
for any finite block-sequence $\tilde{y}=(y_1,\ldots,y_k)$ with
$\max(supp(y_k)) < \min(supp(x_{k+1}))$, the space
$[\tilde{y}^{\frown} (x_n)_{n \geq k+1}]$ is $(k+cC)c^2$-complemented in
$X$: we just
apply the previous fact noting that
$[(x_n)_{n \geq k+1}]$ is $cK$-complemented in $X$ and that the
$k$-dimensional space
$[\tilde{y}]$ is $k$-complemented in $X$.

Now as before we find $K$ such that the set of $K$-complemented block-subspaces
is non-meagre, so residual in some open set $N(\tilde{u})$. By the above,
we have a uniform
control on the constant when we modify a given number of vectors of the
subspace, so
as in Proposition \ref{prop5}, we find $K'$ such that the set of block-sequences
spanning a $K'$-complemented subspace of $X$ is residual in $bb_d(X)$.
Let $\{\tilde{x}_n\}$ be a sequence given by Proposition \ref{prop6}; any
block-sequence
passing through infinitely many $\tilde{x_n}$'s spans a $2K'$-complemented
block-subspace.

Let now $Y=[y]$ be an arbitrary block-subspace of $X$. We easily find
successive intervals of integers $\{E_i, i \in \N\}$, a subsequence
$(\tilde{x}_i)$ of $(\tilde{x}_n)$,
 and a sequence $(\tilde{y}_i)$ of finite
block-sequences  with
 $\tilde{y}_1^\frown \tilde{y}_2^\frown \ldots=y$, such that for all $i$,
$supp(\tilde{x_i}) \cup supp(\tilde{y}_i) \subset E_i$.
By construction, $[\tilde{y}_1^\frown \tilde{x}_{2}^\frown
\tilde{y_3}^\frown \ldots]$
is $2K'$-complemented in $X$, so $Y_1=[\tilde{y}_{1}^{\frown}
\tilde{y}_3^\frown \ldots]$
is $2cK'$-complemented in $X$.
Likewise, $Y_2=[\tilde{y}_{2}^{\frown} \tilde{y}_4^\frown \ldots]$
is $2cK'$-complemented in $X$.
Applying the fact stated at the beginning of the proof, we deduce that
the space $Y=Y_1 \oplus Y_2$ is $4c^3 K'$-complemented in $X$.
\pff

\
Note that we needed to assume unconditionality for our result. If we don't,
we at least
get unconditionality in some block-subspace. Indeed, if every
block-subspace is complemented,
then every block-subspace $Y=[y_n]_{n \in \N}$ has a non-trivial
complemented subspace,
$[y_{2n}]_{n \in \N}$ for example. So $X$ contains no hereditarily decomposable
subspace, and by
Gowers' dichotomy theorem, it must contain an unconditional block-sequence.

It follows from Proposition \ref{uniformity of complementation} and from
the results of
Maurey, Milman and Tomczak-Jaegermann in \cite{MMT} that $c_0$ and $l_p$ are
 the only subsymmetric
unconditional examples.

\begin{coro} \label{asymptotic} Let $X$ be a space with an unconditional basis
$\{e_i\}$
 such that
every block-subspace of $X$ is complemented. Then $X$ is asymptotically
$c_0$ or $l_p$.
\end{coro}

\pf Again $c$ denotes the unconditional constant of $\{e_i\}$. Let $K$ be
given by Proposition \ref{uniformity of complementation}.
 It
follows that all finite block-sequences span $cK$-complemented subspaces of
$X$ (this
could also have been proved directly in the spirit of \cite{DDS}).

With Maurey, Milman, Tomczak-Jaegermann's terminology, in particular
all permissible subspaces of $X$ far
enough are $cK$-complemented. We may then apply their Theorem 5.3 in the
unconditional case
(Remark 5.4.1) to deduce that $X$ is asymptotically $c_0$ or $l_p$.
\pff

\begin{prop}
Let $\{e_i\}$ be an unconditional basic sequence asymptotically $c_0$ or
$l_p$ such that every subsequence of $\{e_i\}$ has a block-sequence
equivalent to $\{e_i\}$ (in particular,
if $\{e_i\}$ is
 subsymmetric or equivalence-block minimal), then
$\{e_i\}$ is
equivalent to the unit basis of $c_0$ or $l_p$.
\end{prop}

\pf
Let $p$ be such that $X=[e_i]_{i \in \N}$ is
asymptotically $l_p$ (the case of $c_0$ is similar).
Assume every subsequence of $\{e_i\}$ has a block-sequence equivalent to
$\{e_i\}$. Then as shown we may (by passing to a subsequence) assume that for
some $C \geq
1$, every subsequence of $\{e_i\}$ has a block-sequence $C$-equivalent to
$\{e_i\}$.

We fix $n \in \N$ and build a winning strategy for Player II in the
asymptotic game of infinite length for producing
a block-sequence $2C$-equivalent to
$(e_i)$. This strategy may then be opposed
to a winning strategy for Player I for producing
length $n$ block-sequences $C'$-equivalent to the unit basis of
$l_p^n$. We get that $e_1,\ldots,e_{n}$
 is $2CC'$-equivalent to $l_p^n$ for all $n$ which will conclude the proof.

Let ${\ku A}=\mgdv (n_k)\in \ram: \exists (x_k)\in bb_N(X)\; \forall k\;
n_{2k}<x_k<n_{2k+1}
\og (x_k) \sim^C (e_k)\mgdh$. We claim that any sequence $(m_k)\in \ram$
contains a further subsequence in $\ku A$, for we can suppose that
$I_k=]m_{2k},m_{2k+1}[\neq\tom$ for all $k$ and therefore take
$(y_i)\sim^C(e_i)$ with $supp(y_i)\subset \for I_n$. Then in between $y_i$
and $y_{i+1}$ there are $n_{2i+1}:=m_{2k+1}$ and  $n_{2i+2}:=m_{2k+2}$,
whereby $(n_k)\in \ku A$.
 So by the infinite Ramsey theorem there is some infinite $A\subset \om$
such that $[A]^\om\subset \ku A$ and there is a $C$-measurable
$f:[A]^\om\til bb_N(X)$ choosing witnesses $(x_k)$ for being in $\ku A$.

Let $\Delta=(\delta_n),\; \delta_n>0$ be such that if two normalized block
sequences are less than $\Delta$ apart, then they are $2$-equivalent. We
choose inductively $n_i<m_i<n_{i+1}$ and sets $B_i\subset\om$ such that
$n_i,m_i\in B_i\subset B_{i-1}/m_{i-1}$ and such that for all $C,D\in
\intv\{n_{j_1},m_{j_1},\ldots,n_{j_k},m_{j_k}\},B_{j_k+1}\inth$
($j_1<\ldots<j_k$)  we have $\|f(C)(k)-f(D)(k)\|<\delta_k$. This can be
done as the unit sphere of $[e_{n_{j_k}+1},\ldots,e_{m_{j_k}-1}]$ is
compact.

 Now in the asymptotic game of infinite length, we can demand that I plays
numbers from the sequence $(n_i)$ and then II replies to
$n_{j_1},\ldots,n_{j_k}$ played by I with some $n_{j_k}<x<m_{j_k}$ such
that for all $C\in
\intv\{n_{j_1},m_{j_1},\ldots,n_{j_k},m_{j_k}\},B_{j_k+1}\inth$ we have
$\|f(C)(k)-x\|<\delta_k$.

Then in the end of the infinite game, supposing that I has played
$(n_{j_k})$ and II has followed the above strategy responding by $(x_k)$,
we have  $n_{j_k}<x_k<m_{j_k}$. Let $(y_k):=f(\{n_{j_k},m_{j_k}\}_\om)$,
then $\|x_k-y_k\|<\delta_k$ and $(x_k)\sim^2(y_k)\sim^C(e_i)$, so II wins.
 \pff

\

\

A space is said to be {\em countably homogeneous} if it has at most
countably many non-isomorphic subspaces. By Corollary \ref{minimality},
a countably homogeneous space
has a block-minimal subspace with an unconditional basis, and one easily
diagonalizes to get
a hereditarily block-minimal subspace.

\begin{prop} \label{prop16} Let $X$ be a coun\-tably homo\-geneous,
here\-di\-ta\-rily block-mini\-mal Banach space with an un\-condi\-tional
basis. Then
ele\-ments in the resi\-dual class of isomorphism for $bb_d(X)$ are iso\-morphic
 to a (possi\-bly infi\-nite) direct sum
of an ele\-ment of each class.
\end{prop}

\pf We write the proof in the denumerable case.
We partition $X$ in a direct sum of subspaces $X_n,n \in \N$ by
partitioning the basis. So each $X_n$
embeds into $X$.
 For each $n$,
choose a representative $E_n$ of the $n$-th isomorphism class
which is a block of $X_n$ (it is possible
because $X_n$ is block-minimal as well). By applications of Gowers' Theorem
  in each $X_n$,
we may pick each vector forming the basis of each $E_n$ far enough, to
ensure that
$E=\sum_{n \in \N}\oplus E_n$ is a block-subspace of $X$.
 We show that $E$ is in the residual class $\cal A$.
 Indeed, if $m$ is such that $E_m \in \cal
A$, then $E \simeq E_m \oplus \sum_{n \neq m} \oplus E_n \simeq E_m$  by
Corollary \ref{coro8}.
 \pff

It follows from the proof above that for any two block-subspaces $A$ and $B$ of
$X$, $A \oplus B$ may be embedded as a block-subspace of $X$; i.e., under the
assumptions of Proposition \ref{prop16},
isomorphism classes of block-subspaces of $X$ form a countable (commutative)
semi-group.

\

We now want to generalize the previous results, considering more general
types of
subspaces. So we fix $X$ to be a Banach space with a basis $\{e_n\}$ and we
denote
by $c$ the  constant of the basis.

 We first notice that we could consider subspaces generated by disjointly
supported
(but not necessarily successive) vectors, and get similar results: we find
a residual class characterized by a "passing through" property.
It follows that
in a non ergodic Banach space with an unconditional basis,
  subspaces generated by disjointly supported vectors
embed complementably  in a given element of the residual class.
We recall that if a Banach space has an unconditional basis such that any
subspace
generated by disjointly supported
vectors is complemented, then the space is $c_0$ or $l_p$ (\cite{LT});
however there does exist
a space $X_U$ cited above from \cite{P}, not isomorphic to $c_0$ or $l_p$, and
 such that
every subspace generated by disjointly supported vectors embeds
complementably in $X_U$.

\

Then we want to represent any subspace of $X$, possibly up to small
perturbations, on
the basis $\{e_n\}$, and get similar results as for the case of block subspaces.
We shall call {\em triangular sequences of blocks} the normalized sequences
of (possibly infinitely supported)
 vectors
in the product
 $X^{\omega}$'s, satisfying for all $k$, $\min(supp(x_k))<\min(supp(x_{k+1}))$
equipped
with the product of the norm topology on $X$. The set of
triangular sequences
 of blocks will be denoted by
$tt$. By Gaussian elimination method, it is clear that any subspace of
 $X$
may be seen as the closed linear space generated by some sequence of $tt$.
Once again it is possible to discretize the problem by considering the set
$tt_d$ of sequences
of vectors in $\bf {Q}$ such
that for all $k$, $\min(supp(x_k))<\min(supp(x_{k+1}))$, and by showing
that for any $x \in tt$ and any $\e>0$, there exists
$x_d$ in $tt_d$ such that $[x_d]$ is $1+\e$-isomorphic to $[x]$.

Our usual method does generalize to this setting. However the
characterization of a residual
set turns out to be only expressed
 in terms of particular subspaces of $tt$, namely those with a finite
dimensional
decomposition on the basis (or F.D.D.). So it gives more information, and
it is actually
easier, to work
 directly with
spaces with a finite dimensional decomposition on the basis.

Note that a space with a F.D.D. on the basis must have the bounded
approximation property.
So our methods only allow us to study subspaces with that property.
In fact, it is easy to check that the set of sequences in $tt_d$ spanning a
space with an F.D.D. on the
basis is residual in $tt_d$, and so, the set of spaces without the
(bounded) approximation property is meagre in our
topology, which explains why with our methods,
we do not seem to be able to "see" spaces without the approximation
property.

  We
 say that two finite-dimensional spaces $F$ and $G$
are {\em successive}, and write $F<G$, if for any $x \in F$, $y \in G$, $x$
and $y$ are successive.
A space with a {\em finite dimensional decomposition on the basis} is a
space of the form
$\oplus_{k \in \N} E_k$,
with successive, finite-dimensional spaces $E_k$; such a space {\em passes
through} $E$ if
$E_k=E$ for some $k$.
We let $fdd$ be the set of infinite sequences of successive
finite-dimensional subspaces, and $fdd_d$ be the set of infinite
sequences of
successive finite-dimensional subspaces which
are spanned by a collection
of vectors with rational
coordinates - equipped with the product of the discrete topology on the set
of finite collections of rational vectors.

\begin{theo} \label{theo19}
Let $X$ be a Banach space with a basis. Then
either $X$ is ergodic, or there exists $K \geq 1$,
 and a sequence
of successive finite dimensional spaces $F_n$ such all spaces with
finite dimensional decomposition on the basis passing
through infinitely many $F_n$'s are mutually $K$-isomorphic.
\end{theo}

\pf Most of the previous proof may be taken word by word; instead of
working with
block-subspaces, i.e. subspaces with $1$-dimensional decomposition on the
basis, we work with
subspaces with arbitrary F.D.D. on the basis; we just have to define the
$\Delta$-expansion of a subspace with a F.D.D. on the basis using an
appropriate distance
between finite-dimensional subspaces. We then end up with characterization
in terms of
"passing through" some finite sequences of finite-dimensional spaces
$\{\tilde{E}_i, i \in \N\}$,
 and this may be simplified to get Theorem \ref{theo19}, noting that we may
choose these sequences to be of length
$1$ (replace each $\tilde{E}_i=(E_i^1,\ldots,E_i^{n_i})$ by $E_i^1 \oplus
\ldots \oplus E_i^{n_i}$).
\pff

\

If the basis is unconditional then we can use the many projections to get
additional results
concerning the residual class.

\begin{prop} \label{prop20}
Let $X$ have an unconditional basis and be non-ergodic. Let $A$ belong to
the residual class in $fdd_d$.
Then every subspace with a finite dimensional decomposition on the basis embeds
 complementably in $A$. Either $A$  fails to have l.u.s.t. or $X$ is $C$-$l_2$
saturated.
 \end{prop}

 \pf  Let $(F_n)$ be given
 by Theorem \ref{theo19}. let $Y=\sum \oplus B_n$ be an arbitrary space
 with a finite dimensional decomposition on the basis.
Passing to a subsequence of $(F_n)$, we may assume that there is
 a partition of $\N$ in successive intervals
$J_i, i \in \N$, such that for all $i$,
$$F_{i-1} < B'_i=\oplus_{n \in J_i} B_n < F_{i+1}.$$
We have that $$A \simeq (\sum \oplus F_i)_{i \in \N} \simeq (\sum \oplus
F_{2k-1})_{k \in \N}
 \oplus (\sum \oplus
F_{2k})_{k \in \N},$$ and so it follows that
$$A \oplus Y \simeq  (\sum \oplus F_{2k-1})_{k \in \N} \oplus (\sum \oplus
B'_{2k})_{k \in \N}
\oplus  (\sum \oplus F_{2k})_{k \in \N} \oplus (\sum \oplus B'_{2k-1})_{k
\in \N}$$
 $$\simeq A \oplus A \simeq A.$$

  For the last part of the Proposition, if $X$ is $C$-$l_2$ saturated for
no $C$, then
 it follows from the Theorem of Komorowski and Tomczak-Jaegermann
(\cite{KT1},\cite{KT2}, and \cite{T} for a
 survey and an improved result) that
 $X$ has a finitely supported subspace $L_n$ without $n$-l.u.s.t.
  for each $n$: indeed Komorowski-Tomczak-Jaegermann's result
 takes care of the finite cotype case, and if $X$
hasn't finite cotype then it contains $l_{\infty}^n$'s uniformly; so it contains
$1+\e$-isomorphic, finitely supported if you wish, copies of any finite
dimensional space, for example without $n$-l.u.s.t..
We may by Theorem \ref{theo19} extend $L_n$ to a space with a
finite dimensional decomposition on the basis
which is $K$ isomorphic to $A$.
Then if $c$ is the unconditional constant of the basis,
$A$ must fail $n/cK$-l.u.s.t. and as $n$ was arbitrary, $A$ fails l.u.s.t..
 \pff

\

Note that it is a consequence of the solution of Gowers and Komorowski-Tomczak to
the Homogeneous Banach
Space Problem, that the following strengthening holds:
a Banach space isomorphic to all its subspaces with the Bounded
Approximation Property
must be isomorphic to $l_2$. This can also be seen as a consequence of the
previous Proposition.
  Note
also that by results of Kadec and Pelczynski (Theorem 2.d.8. in
\cite{LT}), there exists
 a Banach space which is
complementably universal for the set of Banach spaces with the Bounded
Approximation Property; but as mentioned
before, there is  no complementably universal space for the class of
separable Banach spaces. The Bounded
Approximation Property
 seems to draw the
fine line for positive results with our methods.

We conjecture that $l_2$ is the only non-ergodic Banach space. However,
we are not even able to prove that $c_0$ and $l_p,p \neq 2$ are ergodic -
although it is known that those spaces have at least $\omega_1$
non-isomorphic subspaces (\cite{LT}).
 To answer this
question, evidently one would have to consider other types of subspaces
than those generated
by successive blocks,
or disjointly supported blocks; one could consider spaces of the form
$(\sum_{n=1}^{+\infty} \oplus B_n)_p$, with carefully chosen finite-dimensional
$B_n$ so that
this direct sum is isomorphic to a subspace of $l_p$, and play with the
possible choices for
$(B_n)$ (see \cite{LT}, Proposition 2.d.7.).

On the other hand, it is more relevant to restrict the question of
ergodicity to block-subspaces,
 if one
is looking for a significant dichotomy between "regular" and "wild" spaces
with a basis:
 in this setting, $c_0$ and
$l_p$ are, as they should be, on the regular (i.e. non-ergodic) side of the
dichotomy.

\paragraph{Acknowledgements}
We thank A. Pelczar and N. Tomczak- Jaegermann for drawing our attention to
the Theorem
of Maurey, Milman and Tomczak- Jagermann in \cite{MMT} which led us to
Proposition \ref{uniformity of complementation}
and Corollary \ref{asymptotic}. We also thank S. Todorcevic for a useful
discussion about asymptotic games.


\begin{thebibliography}{AAA}



\bibitem{B} B. BOSSARD, {\em A coding of separable Banach spaces.
Analytic and coanalytic families of Banach
spaces}, Fundamenta Mathematicae, 172 (2002), {\bf 2}, 117-152.

\bibitem{BL} J. BAGARIA, J. LOPEZ-ABAD, {\em Weakly Ramsey Sets in
Banach Spaces}, Advances in Mathematics 160,
{\bf 2} (2001), 133-174.

\bibitem{CS} P.G. CASAZZA, T. SHURA, {\em Tsirelson's space, Lecture
Notes in Mathematics, 1363, Springer Verlag, (1989).}

\bibitem{DDS} W.J. DAVIS, D.W. DEAN, I. SINGER, {\em Complemented
subspaces and $\Lambda$-systems in Banach spaces}, Israel Journal of Mathematics
{\bf 6} (1968), 303-309.




\bibitem{FR} V. FERENCZI, C. ROSENDAL: {\em On the number of non
isomorphic subspaces of a Banach space},
preprint.

\bibitem{G0} W.T. GOWERS, {\em A solution to Banach's hyperplane
problem}, Bulletin of the London
Mathematical Society {\bf 26} (1994), 523-530.

\bibitem{G1} W.T. GOWERS, {\em A new dichotomy for Banach spaces},
 Geometric and Functional Analysis {\bf 6} (1996), 6, 1083-1093.



\bibitem{G2} W.T. GOWERS, {\em An infinite Ramsey Theorem and some
Banach-space
dichotomies}, Annals of Math. 156 (2002), 3.



\bibitem{GM} W.T. GOWERS and B. MAUREY, {\em The unconditional basic
sequence problem}, J.Amer.Math.Soc. {\bf 6} (1993), 851-874.

\bibitem{GM2} W.T. GOWERS and B. MAUREY, {\em Banach spaces with small
spaces
of operators}, Math. Ann. {\bf 307} (1997), 4, 543-568.

\bibitem{just} W. JUST and M. WEESE, {\em Discovering modern set theory. II, Set-theoretic tools for every mathematician}, American mathematical Society (1997).


\bibitem{kal} N. KALTON, {\em A remark on Banach spaces isomorphic to their squares}, Contemporary  Mathematics, Volume {\bf 232} (1999).



\bibitem{kec} A. KECHRIS, {\em Classical descriptive set theory},
Sprin\-ger-Ver\-lag,
 New York (1995).



\bibitem{KT1} R. KOMOROWSKI and N. TOMCZAK-JAEGERMANN, {\em Banach spaces
without local unconditional structure}, Israel Journal of Mathematics {\bf 89}
(1995), 205-226.



\bibitem{KT2} R. KOMOROWSKI and N. TOMCZAK-JAEGERMANN, {\em Erratum to:
"Banach spaces without local unconditional structure"},
 Israel Journal of Mathematics {\bf 105}
(1998), 85-92.



\bibitem{LT} J. LIN\-DEN\-STRAUSS and L. TZA\-FRI\-RI,
{\em Clas\-si\-cal Ba\-nach spa\-ces I}, Sprin\-ger-Ver\-lag,
 New York (1977).

\bibitem{JLA} J. LOPEZ-ABAD,
{\em Coding into Ramsey sets}, preprint.

\bibitem{MMT} B. MAUREY, V.D. MILMAN and  N. TOMCZAK-JAEGERMANN,
{\em Asymptotic infinite-dimensional theory of Banach spaces},
Geometric Aspects of Functional Analysis, Israel, 1992-1994, 149-175,
Oper. Th. Adv. Appl., Birkhauser, Basel 1995.

\bibitem{OS} E. ODELL and T. SCHLUMPRECHT, {\em The distortion problem},
Acta Math. {\bf 173} (1994), no. 2, 259-281.

\bibitem{P} A. PELCZYNSKI, {\em Universal bases}, Studia Mathematica,
{\bf 40} (1971), 239-242.



\bibitem{R} C. ROSENDAL: {\em On the complexity of isomorphism of
separable Banach spaces}, preprint.



\bibitem{S} A. SZANKOWSKI, {\em Subspaces without the Approximation
Property}, Israel Journal of Mathematics {\bf 30} 1-2 (1978), 123-129.



\bibitem{T} N. TOMCZAK-JAEGERMANN, {\em Banach spaces with many
  isomorphisms}, Semin\'ario Brasileiro de An\'alise {\bf 48}, Petropolis,
  Rio de Janeiro,  November 1998, 189-210.



\end{thebibliography}
\end{document}